\documentclass{jnmp01b}

\usepackage{amsmath}
\usepackage{graphicx}

\numberwithin{equation}{section}

\newtheorem{prop}{Proposition}

\begin{document}

\renewcommand{\evenhead}{B.A.\ Springborn}
\renewcommand{\oddhead}{The Toy Top, an Integrable System 
of Rigid Body Dynamics}


\thispagestyle{empty}

\begin{flushleft}
\footnotesize \sf
Journal of Nonlinear Mathematical Physics \qquad 2000, V.7, N~3,
\pageref{firstpage}--\pageref{lastpage}.
\hfill {\sc Article}
\end{flushleft}

\vspace{-5mm}

\copyrightnote{2000}{B.A.\ Springborn}

\Name{The Toy Top, an Integrable System of Rigid Body Dynamics}

\label{firstpage}

\Author{Boris A.\ SPRINGBORN}

\Adress{Technische Universit\"at Berlin,
Fachbereich Mathematik, Sekr. MA 8-5,\\
Strasse des 17. Juni 136,
D-10623 Berlin, Germany \\
E-mail: springb@sfb288.math.tu-berlin.de}

\Date{Received March 20, 2000; 
Accepted May 2, 2000}

\begin{abstract}
\noindent
A toy top is defined as a rotationally symmetric body moving
in a constant gravitational field while one point on the symmetry axis
is constrained to stay in a horizontal plane. It is an integrable
system similar to the Lagrange top. 
Euler-Poisson equations are derived. Following Felix Klein, the special
unitary group $\rm SU(2)$ is used as configuration space and the
solution is given in terms of hyperelliptic integrals. The curve
traced by the point moving in the horizontal plane is analyzed, and a
qualitative classification is achieved. The cases in which the
hyperelliptic integrals degenerate to elliptic ones are found and the
corresponding solutions are given in terms of Weierstrass elliptic
functions.  
\end{abstract}



\section{Introduction}

The three famous integrable cases of rigid
body motion, the tops of Euler, Lagrange and Kowalewsky,
have been paramount examples in the theory of integrable
systems. The modern algebra-geometric approach, using
Lax pairs with a spectral parameter \cite{IntSysI},
\cite{IntSysII} has been applied to all three. 
It is surprising
that the following system appears only sporadically in the
literature: A rotationally symmetric rigid body moving in a
homogeneous gravitational field with one point on its axis not fixed,
but constrained to move in a horizontal plane. Following F. Klein
\cite[p.\ 58]{Klein:Lectures}, we call such a system a {\em toy
top}.

In many ways, the toy top is similar to Lagrange's top. It is
completely integrable due to the same kind of rotational symmetry. 
The solution leads to hyperelliptic integrals instead of
elliptic ones, but their analytic properties are similar.
It seems that Poisson is the first to solve the
system \cite{Poisson:Mechanik}, using Euler angles. It is treated
similarly by E. T. Whittaker \cite{Whittaker} and F. Klein
\cite{Klein:Mechanik}. Later, Klein discovered that, as in the case of
Lagrange's top, simpler solutions are obtained when $\rm SU(2)$ is used as
configuration space \cite{Klein:Buch}, \cite{Klein:Lectures}.

In this paper we close some gaps in the classical
treatment. Following A. I. Bobenko's and Yu.\ B. Suris' treatment of
the Lagrange top \cite{BobenkoSuris1}, \cite{BobenkoSuris2}, new
equations of
motion are derived within the framework of Lagrangian mechanics on Lie
groups. 

As the toy top moves, its tip traces a curve on the supporting
plane. Formulas for this curve are derived and qualitatively different
cases are classified. 

For certain values of the first integrals, the hyperelliptic
integrals appearing in the solution degenerate to elliptic
integrals. These cases are classified and the corresponding solutions
are given in terms of Weierstrass elliptic functions. 


\section{Preliminaries}

After it is shown how the group of rotations can be considered
the configuration space of the toy top, the alternative use of $\rm
SU(2)$ is discussed. Finally, M\"obius transformations that rotate the
Riemann sphere of numbers are discussed for use in chapter \ref{sec:tip}. 


\subsection{The group $\rm SU(2)$ as configuration space}

The configuration space of the toy top is the
direct product of the rotation group with the group of translations in
the plane. But since both the gravitational and the resistive force
of the plane are vertical, the horizontal component
of the velocity of the top's center of mass is constant. After 
changing to a
suitable moving coordinate system if necessary, 
the center of mass moves only vertically. This reduces the
configuration space to the rotation group: The top can be brought into
any admissible position by a rotation around its center of mass,
followed by a vertical translation. But the latter is determined by the
former due to the constraint. 

Instead of the special orthogonal group $\rm SO(3)$ or Euler angles, the
special unitary group $\rm SU(2)$ will be used to describe
configurations of the toy top. As F. Klein discovered, this leads to
simpler formulas.
The special unitary group of two dimensions is the matrix group
\begin{equation*}
{\rm SU(2)}=\left\{\left. 
\begin{pmatrix} \alpha & \beta \\ \gamma & \delta \end{pmatrix} 
\in\mathbb{C}^{2\times 2}\;\right|\;
\alpha\delta-\beta\gamma=1,\;\delta=\bar{\alpha},\;\gamma=-\bar{\beta}
\right\}.
\end{equation*}
Its Lie algebra consists of the skew hermitian matrices with trace
zero:
\begin{equation*}
\mathrm{su}(2)=\left\{\left. 
\begin{pmatrix} ix_3 & -x_2+ix_1 \\ x_2+ix_1 & -ix_3 \end{pmatrix}
\in\mathbb{C}^{2\times 2}\;\right|\; 
x_1, x_2, x_3 \in \mathbb{R}
\right\}.
\end{equation*}
Identify $\mathbb{R}^3$ with $\mathrm{su}(2)$ via
\begin{equation*}
(x_1,x_2,x_3)\mapsto \frac{1}{2}
\begin{pmatrix} ix_3 & -x_2+ix_1 \\ x_2+ix_1 & -ix_3 \end{pmatrix}.
\end{equation*}
This corresponds to the choice of basis
\begin{equation*}
e_1=\frac{1}{2}\begin{pmatrix} 0 & i \\ i & 0 \end{pmatrix},\quad
e_2=\frac{1}{2}\begin{pmatrix} 0 & -1 \\ 1 & 0 \end{pmatrix},\quad
e_3=\frac{1}{2}\begin{pmatrix} i & 0 \\ 0 & -i \end{pmatrix}.
\end{equation*}
The factor $1/2$ is introduced so that the cross product in
$\mathbb{R}^3$ corresponds to the Lie bracket in $\mathrm{su}(2)$:
\begin{equation*}
x\times y=[x,y].
\end{equation*}
The induced scalar product on $\mathrm{su}(2)$ is
\begin{equation*}
\langle x,y \rangle =-2\,\mathrm{Tr}(x\,y).
\end{equation*}
The adjoint action of $\mathrm{SU}(2)$ on its Lie algebra is
orthogonal and orientation preserving. The homomorphism 
\begin{align*}
\mathrm{Ad}:{\mathrm{SU} (2)} &\rightarrow\mathrm{SO}(3) \\
    \Phi&\mapsto (\mathop{\mathrm{Ad}}\nolimits_\Phi: x\mapsto\Phi
x \Phi^{-1})
\end{align*}
is 2 to 1 with Kernel $\{\pm 1\}$.

Consider two orthonormal coordinate systems whose origin is the center
of mass of the toy top. The first, called the fixed frame, has axes
whose directions are fixed. Assume that the third axis points
upward. The second coordinate system, called the moving frame, moves
with the toy top. Assume its third axis points along the top's
symmetry axis away from the top's tip. The momentary position of the
top is given by a matrix $\Phi\in{\rm SU}(2)$
such that if $X=(X_1, X_2, X_3)$ is the coordinate vector of some
point in the body frame then its coordinates $x=(x_1, x_2, x_3)$ with
respect to the fixed frame are given by 
\begin{equation}\frac{1}{2}
\begin{pmatrix} ix_3 & -x_2+ix_1 \\ x_2+ix_1 & -ix_3 \end{pmatrix}
=\frac{1}{2}\Phi
\begin{pmatrix} iX_3 & -X_2+iX_1 \\ X_2+iX_1 & -iX_3 \end{pmatrix}
\Phi^{-1}.
\end{equation} 
In the context of rigid body mechanics, the entries of an
${\rm SU}(2)$ matrix are called {\em Cayley-Klein parameters}. 

Suppose that a curve in ${\rm SU}(2)$ describing the motion of the top 
is given by
\begin{equation*}
\Phi(t)=
\begin{pmatrix} \alpha(t) & \beta(t) \\ \gamma(t) & \delta(t) \end{pmatrix}.
\end{equation*}
A moving point whose coordinate vector is $X(t)$ in the moving frame
has coordinates $x(t)=\Phi(t)X(t)\Phi^{-1}(t)$ in the fixed
frame. Taking the derivative, one obtains 
$x'=\Phi X' \Phi^{-1} + [\omega, x]$
and
$X'=\Phi^{-1} x' \Phi - [\Omega, X]$,
where $\omega=\Phi'\Phi^{-1}$ and $\Omega=\Phi^{-1}\Phi'$. 
Since
$\omega=\Phi\Omega\Phi^{-1}$, they represent the same vector, the {\em
angular velocity vector}, in the fixed and moving frame, respectively.
One obtains
\begin{equation}\label{eqn:omega}
\omega=\frac{1}{2}
\begin{pmatrix}
i\omega_3 & \omega_2+i\omega_1  \\ -\omega_2+i\omega_1  & i\omega_3
\end{pmatrix}=
\begin{pmatrix}
\alpha'\delta-\beta'\gamma & -\alpha'\beta+\beta'\alpha \\
\gamma'\delta-\delta'\gamma & -\gamma'\beta+\delta'\alpha
\end{pmatrix} 
\end{equation}
and
\begin{equation}\label{eqn:Omega}
\Omega=\frac{1}{2}
\begin{pmatrix}
i\Omega_3 & \Omega_2+i\Omega_1  \\ -\Omega_2+i\Omega_1  & i\Omega_3
\end{pmatrix}=
\begin{pmatrix}
 \delta\alpha'-\beta\gamma' &  \delta\beta'-\beta\delta' \\
-\gamma\alpha'+\alpha\gamma' & -\gamma\beta'+\alpha\delta'
\end{pmatrix}.
\end{equation}


\subsection{M\"obius transformations}

There is another way to establish the $\rm SO(3)$ action of $\rm
SU(2)$ via M\"obius transformations of the Riemann sphere. It is
summarized here for reference in Section \ref{sec:tip}. For more
detail see \cite[pp. 29ff]{Klein:Ikosaeder}. 

The unit sphere
$S^2=\{(\xi,\eta,\zeta)\in\mathbb{R}^3\,|\,\xi^2+\eta^2+\zeta^2=1\}$ 
is mapped conformally onto the extended complex plane
$\mathbb{C}\,\cup\{\infty\}$ by stereographical
projection from the north pole. To the point $z=x+iy$ in
the complex plane corresponds the
point in the sphere with coordinates
\begin{equation}\label{eqn:inerse-stereo-proj}
\xi=\frac{2x}{|z|^2+1},\quad \eta=\frac{2y}{|z|^2+1},\quad
\zeta=\frac{|z|^2-1}{|z|^2+1}\;. 
\end{equation}
The conformal
maps of the Riemann sphere onto itself are described by M\"obius
transformations
\begin{equation}\label{eqn:Moeb}
z\mapsto\frac{\alpha z + \beta}{\gamma z +\delta}\,,\quad 
\text{where } \alpha\delta-\beta\gamma=1.
\end{equation}
The map sending a matrix in SL(2,$\mathbb{C}$) with entries $\alpha$, $\beta$,
$\gamma$, $\delta$ to the M\"obius transformation (\ref{eqn:Moeb})
is a group homomorphism with kernel $\{1,-1\}$.
In particular,
the isometric transformations of the sphere, i.e.\ the rotations, are
conformal. They correspond to those M\"obius transformations with
$\delta=\bar{\alpha}, \gamma=-\bar{\beta}$. This is the image of 
${\rm SU}(2)$ under the above homomorphism.


\section{The Lagrangian, equations of motion, and their solution in
terms of hyperelliptic functions} 
\label{sec:sysandsol}

The Lagrangian description of the toy top is used to derive
equations of motion. The integrals of motion---total energy and two
further integrals connected to the rotational symmetries of the
system---are used
to reduce the system to one degree of freedom. It is then solved in
hyperelliptic integrals. Most of the results in this section, but not
the equations of motion (\ref{eqn:motion}) and their derivation, are
already contained in the classical works mentioned in the
introduction. 


\subsection{The Lagrangian, first integrals, and reduced equation of
motion} 

In general, the Lagrangian is a function on the tangent bundle of the
configuration space. If this space is a Lie group, it is
convenient to trivialize the tangent bundle via left or right
multiplication. We will use right multiplication:
\begin{align*}
\mathrm{SU}(2)\times\mathrm{su}(2) &\rightarrow\mathrm{TSU}(2) \\
(\Phi,\omega)&\mapsto (\Phi,\omega\Phi).
\end{align*}
This corresponds to a description in the fixed
frame. For a general exposition of this approach to mechanical systems
similar to the Lagrange top, see \cite{BobenkoSuris1},
\cite{BobenkoSuris2}.

In these fixed frame coordinates, the kinetic and potential energy
functions are
\begin{equation*}\label{eqn:kinetic}
T(\Phi,\omega)=
\frac{1}{2} A     \langle \omega,\omega \rangle +
\frac{1}{2} (C-A) \langle \omega,r \rangle^2 +
\frac{1}{2}\, ps \langle [r,k],\omega \rangle^2,
\end{equation*}
and
\begin{equation*}
V(\Phi,\omega)= p \langle r,k \rangle,
\end{equation*}
where $k=e_3$ is the unit vector pointing
vertically upward, $r=\Phi k\Phi^{-1}$ is the unit vector
pointing in the direction of the top's axis, $s$ is the
distance between the tip of the top and its center of mass, $p$ is
the product of $s$ and the mass, and $C$ and $A$ are the inertia
moments of the top with respect to 
the symmetry axis and any perpendicular axis through the center of
mass. The
gravitational acceleration is assumed to be one, which can be achieved
by a suitable choice of units, e.g.\ for mass. Note that 
$r'=[\omega,r]$. The kinetic energy of the toy top is composed of
two parts: one (the first two terms in equation (\ref{eqn:kinetic})) is
due to the rotation of the top around its center of mass and the other
is due to the vertical motion of the center of mass. If this last term
$(ps/2) \langle [r,k],\omega \rangle^2$ was missing, and the inertia
moments were related to the tip and not the center of mass, one would
obtain the kinetic energy of Lagrange's top.
The Lagrange function is 
\begin{equation*}\label{eqn:lagrange}
L=T-V,
\end{equation*} 
with corresponding momentum
\begin{equation*}
m=\frac{\partial L}{\partial \omega}=
A\omega +
(C-A)\langle \omega,r \rangle r +  
ps \langle [r,k],\omega \rangle [r,k].
\end{equation*}
It is convenient to introduce a derivative
$DL(\Phi,\omega)\in\mathrm{su}(2)$ with respect to the first argument
of $L$ by
\begin{equation*}
\langle DL(\Phi,\omega),\eta \rangle =
\frac{d}{d\epsilon}
L(e^{\epsilon\eta}\Phi,\omega)\big|_{\epsilon=0}.
\end{equation*}
One obtains
\begin{equation*}
DL(\Phi,\omega)=(C-A)\langle \omega,r \rangle [r,\omega] +
ps\langle [r,k],\omega \rangle [r,[k,\omega]] + p[r,k].
\end{equation*}

\begin{prop}
The equations of motion are
\begin{equation}\label{eqn:motion}
\begin{split}
m' &= [\omega,m] + DL(\Phi,\omega) \\
\Phi' &= \omega\Phi.
\end{split}
\end{equation}
The first integrals
\begin{equation*}
h=T+V,\quad l=\langle m,k \rangle,\quad\text{and}\quad n=\langle m,r \rangle
\end{equation*}
allow the reduction to 
\begin{equation}\label{eqn:reduced}
\frac{s}{2}\left(\frac{A}{ps}+1-u^2\right){u'}^2 =
\frac{1}{p}\left(1-u^2\right)\left(h-\frac{1}{2C}n^2-pu\right) -
\frac{1}{2Ap}(l-nu)^2,
\end{equation}
where $u=\langle r,k \rangle$.
\end{prop}
\begin{proof}
Consider a variation $\Phi(t,\epsilon)$ of $\Phi(t,0)$ with fixed
end-points $t=t_0$, $t=t_1$. Suppose 
$\partial\Phi/\partial t=\omega\Phi$ (yielding the second equation of
motion) and 
$\partial\Phi/\partial \epsilon=\eta\Phi$.
For the variation of the action functional
$S=\int_{t_0}^{t_1} L\,dt$ one obtains,
\begin{equation*}
\left.\frac{dS}{d\epsilon}\right|_{\epsilon=0}=
\int_{t_0}^{t_1}\langle
DL(\Phi,\omega)-[m,\omega]-m',\eta\rangle\,dt.
\end{equation*}
(Note that 
$\partial\omega/\partial\epsilon-\partial\eta/\partial
t+[\omega,\eta]=0$.) This yields the first equation of motion. The
total energy $h$ is always a first integral; $l$ and $n$ are easily
checked by direct calculation. The integral $l$ is due to the
rotational symmetry around the vertical axis and $n$ is due to the
rotational symmetry around the top's axis.

To achieve the reduction, note that in terms of the angular velocity,
the momentum integrals are
$l=A\langle \omega,k \rangle + (C-A)\langle \omega,r \rangle \langle
r,k \rangle$
and $n=C\langle \omega,r \rangle$.
Thus one obtains
$\langle \omega,r \rangle = n/C$ and
$\langle \omega,k \rangle = l/A +
\left( 1/A-1/C \right)nu$.
Furthermore,
$\langle \omega,[r,k] \rangle=u'$.
Unless $r=\pm k$, one finds the following representation of $\omega$
in the basis $(r,k,[r,k])$: 
\begin{equation*}
(1-u^2)\,\omega = 
\left( C^{-1}(1-u^2)n - A^{-1}(l-nu)u \right) r + A^{-1}(l-nu)k +
u'[r,k].
\end{equation*}
Now one can express the energy integral $h$ in terms of $l$,$n$,$u$
and $u'$. This yields the reduced equation of motion
(\ref{eqn:reduced}). 
\end{proof}


\subsection{New dynamical constants}
\label{sec:dynpar}

We introduce new dynamical constants $e_1$, $e_2$, $e_3$ to replace
$n$, $h$ and $l$. Together with $\pm e_4$, also introduced in this
section, they are the branchpoints of the hyperelliptic Riemann
surface defined in the next section. 
\begin{prop}
The third degree polynomial on the right hand side of
equation (\ref{eqn:reduced}) has three real roots. Two of them lie in
the closed interval $[-1,1]$ and the third is greater or equal to
$1$. Hence, equation (\ref{eqn:reduced}) can be written
\begin{equation}\label{eqn:reduced2}
-\frac{s}{2}\left(u^2-e_4^2\right){u'}^2 =
\left(u-e_1\right)\left(u-e_2\right)\left(u-e_3\right),
\end{equation}
where 
\begin{equation}\label{eqn:eineq}
-1\leq e_1\leq e_2\leq 1\leq e_3
\end{equation}
and 
\begin{equation}\label{eqn:e4}
e_4=\sqrt{1+\frac{A}{ps}}.
\end{equation}
\end{prop}
\begin{proof}
This follows from the fact that for the dynamical constants $l$, $n$
and $h$ to be
physically feasible, there has to be a value for $u$ between $1$
and $-1$ such that equation (\ref{eqn:reduced}) leads to real
$u'$. This means that the polynomial on the right hand side must be
nonnegative somewhere in the interval $[-1,1]$. But at $u=\pm 1$ it
takes the nonpositive values $-(l\mp n)^2/(2Ap)$, and it goes to $\pm
\infty$ for $u\rightarrow \pm\infty$. Hence there have to be three
real zeroes, situated as stated. 
\end{proof}

Conversely, it can be shown that inequality (\ref{eqn:eineq}) is the
only constraint for the zeroes. I.e., given three numbers
satisfying this inequality, there is a state of motion of the top with
$e_1\leq u\leq e_2$ that leads to exactly these zeroes.

Note that the initial value problem with differential equation
(\ref{eqn:reduced2}) and an initial value for $u$ between $e_1$ and
$e_2$ has more than one solution. First, the initial sign of $u'$ has
to be determined. But even then the solution is only unique up to the
moment when $u$ reaches $e_1$ or $e_2$. After that, $u$ can either
reverse its path straight away, or stay constant for some time. Of
course, this indeterminacy is not inherent in the original system, but
is introduced by the reduction. Upon closer examination of the system,
one finds that $u$ is constantly equal to $e$ only if $e$ is a double
zero of the right hand side of
(\ref{eqn:reduced2}). See \cite[pp.\ 280ff]{Klein:Buch} for
a discussion of the analogous situation in the case of the Lagrange
top. Except for this
singular case, which will be treated in Sections \ref{sec:reg-prec}
and \ref{sec:aper}, only those solutions of (\ref{eqn:reduced2}) have
to be considered that oscillate periodically between the isolated
minima and maxima $e_1$ and $e_2$.

From now on, we will use $e_1$, $e_2$ and $e_3$ as dynamical constants
instead of $l$, $n$ and $h$. However, the first do not
uniquely determine the latter. Indeed, substituting $1,-1$ and $0$
into the right hand sides of (\ref{eqn:reduced}) and
(\ref{eqn:reduced2}), one obtains the relations
\begin{equation}\label{eqn:constants} 
\begin{split}
&\left(1-e_1\right)\left(1-e_2\right)\left(1-e_3\right) =
-\frac{1}{2Ap}(l-n)^2 \\
&\left(-1-e_1\right)\left(-1-e_2\right)\left(-1-e_3\right) 
= -\frac{1}{2Ap}(l+n)^2 \\
&-e_1 e_2 e_3 = \frac{1}{p}\left(h-\frac{n^2}{2C}-\frac{l^2}{2A}\right).
\end{split}
\end{equation}
If the constants $e_1, e_2, e_3$ are given and none of them is equal
to $1$ or $-1$, then the first two equations of (\ref{eqn:constants})
have four solutions for $(l,n)$. If $(a,b)$ is one of them, the others
are $(-a,-b), (b,a)$ and $(-b,-a)$. Note that replacing $l$ and $n$ by
$-l$ and $-n$ is equivalent to looking at a mirror image of the
system. Once a solution for $(l,n)$ is
chosen, the third equation of (\ref{eqn:constants}) determines the
value of $h$. Unless $A=C$, the four solutions for $(l,n)$ lead to two
different values for $h$. 

If one of the constants $e_1, e_2, e_3$ is equal to $1$ or $-1$,
then $l=n$ or $l=-n$, respectively, so that there are only two
solutions. 


\subsection{The hyperelliptic time integral}
\label{sec:t-int}

Solving equation (\ref{eqn:reduced2}) by separation of variables, one
obtains for $t$ the hyperelliptic integral 
\begin{equation}\label{eqn:tint0}
t=\int{\sqrt{\frac{-s\left(u^2-e_4^2\right)}
   {2\left(u-e_1\right)\left(u-e_2\right)\left(u-e_3\right)}}\;du}.
\end{equation}
Allow complex values for $u$
and $t$ so that it becomes an Abelian integral on a hyperelliptic
Riemann surface.
\begin{figure}
\begin{center}
\includegraphics{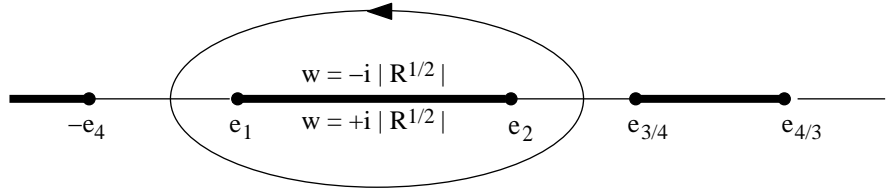}
\end{center}
\caption{The top sheet of the hyperelliptic curve $C$.}
\label{fig:curve1}
\end{figure}
Let $C$ be the hyperelliptic curve given by $w^2=R(u)$, where
\begin{equation*}
 R(u)=(u-e_1)(u-e_2)(u-e_3)(u^2-e_4^2).
\end{equation*}
For future reference we remark that by comparison with equations
(\ref{eqn:reduced}) and (\ref{eqn:reduced2}),
\begin{equation}\label{eqn:R(u)-with-l-n-h}
R(u)=\left(u^2-1-\frac{A}{ps}\right) 
\left(
\frac{1}{p}(1-u^2)\left(h-\frac{1}{2C}\,n^2-pu\right)
-\frac{1}{2Ap}(l-nu)^2
\right).
\end{equation}
The hyperelliptic curve $C$ is a two sheeted branched
covering of the $u$-plane with the six branch points $e_1$, $e_2$,
$e_3$,
$e_4$, $-e_4$, and $\infty$. Figure \ref{fig:curve1} shows one sheet
of $C$ with cuts along the real axis. This will be called the top
sheet. Lift the path
along which $u$ moves during the motion of the top to a path on the
hyperelliptic curve which is homologous to the path drawn in the
figure. Then
\begin{equation}\label{eqn:t}
t=-i\sqrt{\frac{s}{2}}\int\left(u^2-{e_4}^2\right)\frac{du}{w}\,.
\end{equation}
This is an Abelian integral of the second kind which has a simple pole
at infinity: Introducing a holomorphic parameter $v^2=1/u$ around
$\infty$, one obtains the asymptotic expansion
$t=-i\sqrt{2s}v^{-1}+O(1)$ for $v\rightarrow 0$.
 

\subsection{Hyperelliptic integrals for the Cayley-Klein parameters}
\label{sec:abgd-ints}

In the following proposition, hyperelliptic integrals for the
Cayley-Klein parameters are presented and their analytical properties
discussed. 
These solutions are similar to the corresponding results for
the Lagrange top. In the latter case, one obtains elliptic integrals
and not hyperelliptic ones as here, but they have the same kind of
singularities at corresponding places. Also, the reduction to the
case $A=C$ of spherical tops works for the Lagrange top as
well. According to
F. Klein \cite[p. 234]{Klein:Buch}, the reduction of Lagrange's top to
the spherical case  was first noticed by Darboux.

\begin{prop}\label{prop:abgdints}
Suppose first that $A=C$. Then the solution
of the system is given by the hyperelliptic integrals
\begin{equation}\label{eqn:abgdints}
\begin{split}
\log \alpha &= \int \frac{1}{2(u+1)}\left(w+
   \frac{u^2-{e_4}^2}{1-{e_4}^2}w_{-1}\right)\frac{du}{w}\\
\log \beta &= \int \frac{1}{2(u-1)}\left(w-
   \frac{u^2-{e_4}^2}{1-{e_4}^2}w_{+1}\right)\frac{du}{w}\\
\log \gamma &= \int \frac{1}{2(u-1)}\left(w+
   \frac{u^2-{e_4}^2}{1-{e_4}^2}w_{+1}\right)\frac{du}{w}\\
\log \delta &= \int \frac{1}{2(u+1)}\left(w-
   \frac{u^2-{e_4}^2}{1-{e_4}^2}w_{-1}\right)\frac{du}{w}\,.
\end{split}
\end{equation}
The constants $w_{\pm 1}$ denote one of the two values of $w$ on $C$
over $u=\pm 1$, i.e.\
\begin{equation*}
w_{+1}=\pm\sqrt{R(1)}\quad\text{and}\quad w_{-1}=\pm\sqrt{R(-1)}\,.
\end{equation*}
(There are four possible ways to choose, in
accordance with the indeterminacy of $l$ and $n$ in terms of
$e_1,e_2,e_3$; see Section \ref{sec:dynpar}). 

They are Abelian integrals of the third kind, the
differentials under the integral sign having two simple poles
each, with residues $\pm 1$ at the places shown in Table
\ref{tab:anaprop}, and no other singularities.

\begin{table}
\caption{The poles of the logarithmic differentials.}
\label{tab:anaprop}
\small
\begin{center}
\begin{tabular}{c|l|c}
Differential& pole with res. 1 at  & pole with res. -1 at \\
\hline
$d\alpha/\alpha$  & $(u,w)=(-1,w_{-1})$  & $\infty $            \\
$d\beta/\beta$  & $(u,w)=(+1,-w_{+1})$ & $\infty $            \\
$d\gamma/\gamma$  & $(u,w)=(+1,w_{+1})$  & $\infty $            \\
$d\delta/\delta$  & $(u,w)=(-1,-w_{-1})$ & $\infty $            \\
\end{tabular}
\end{center}
\end{table}

If $A\neq C$, the solution is
\begin{equation}
      \Phi \cdot
         \begin{pmatrix}
            \exp( \frac{i}{2} \tau t) & 0 \\
            0 & \exp(-\frac{i}{2} \tau t)
         \end{pmatrix},
\end{equation}
where 
\[\tau=\left(\frac{1}{C}-\frac{1}{A}\right)n.
\]
I.e., the solution differs from the solution for the spherical top
only by a rotation with constant speed around the top's axis.
\end{prop}
\begin{proof}
First, the expressions for the logarithmic differentials of
$\alpha,\beta,\gamma,\delta$ are derived in the general case. The reduction to the
case of spherical tops is then immediate. Finally, the analytic
properties of the differentials are examined. 

Observe that 
\begin{equation}
\begin{split}
1 &= \alpha\delta-\beta\gamma \\
u &= \alpha\delta+\beta\gamma.
\end{split}
\end{equation}
This implies
\begin{equation}\label{eqn:up1um1}
u+1=2\alpha\delta \text{\quad and\quad} u-1=2\beta\gamma,
\end{equation}
such that
\begin{equation}\label{eqn:du}
\begin{split}
du&=2(\delta\,d\alpha+\alpha\,d\delta)\\
du&=2(\gamma\,d\beta+\beta\,d\gamma).
\end{split}
\end{equation}
By equations (\ref{eqn:omega}) and (\ref{eqn:Omega}), the components
of the angular velocity vector in the direction of the vertical axis
and the top's symmetry axis are given by
\begin{align*}
i\omega_3 &= 2(\alpha'\delta-\beta'\gamma)=2(-\alpha\delta'+\beta\gamma'), \\
i\Omega_3 &= 2(\delta\alpha'-\beta\gamma')=2(-\alpha\delta'+\beta'\gamma).
\end{align*}
This implies 
\begin{align*}
i(\Omega_3+\omega_3)dt &=2(\delta\,d\alpha-\alpha\,d\delta), \\
i(\Omega_3-\omega_3)dt &=2(\gamma\,d\beta-\beta\,d\gamma).
\end{align*}
With equations (\ref{eqn:up1um1}) and (\ref{eqn:du}) one obtains
\begin{equation}\label{eqn:du-plus-iOmegas-dt}
\begin{split}
du + i(\Omega_3+\omega_3)dt &= 4\,\delta\,d\alpha = 2(u+1)\frac{d\alpha}{\alpha}\\
du + i(\Omega_3-\omega_3)dt &= 4\,\gamma\,d\beta = 2(u-1)\frac{d\beta}{\beta}\\
du - i(\Omega_3-\omega_3)dt &= 4\,\beta\,d\gamma = 2(u-1)\frac{d\gamma}{\gamma}\\
du - i(\Omega_3+\omega_3)dt &= 4\,\alpha\,d\delta = 2(u+1)\frac{d\delta}{\delta}.
\end{split}
\end{equation}
Now it follows from $n=C\Omega_3$ and $l=A\omega_3+(C-A)\Omega_3 u$ that
\begin{align}\label{eqn:omegasum}
\Omega_3+\omega_3 &= 
\frac{1}{A}(l+n)+(u+1)\left(\frac{1}{C}-\frac{1}{A}\right)n, \\
\Omega_3-\omega_3 &= 
-\frac{1}{A}(l-n)-(u-1)\left(\frac{1}{C}-\frac{1}{A}\right)n.
\end{align}
Note that from equation (\ref{eqn:R(u)-with-l-n-h}), 
\begin{equation*}
R(\pm 1) = \frac{1}{2p^2s}(l\mp n)^2.
\end{equation*}
Hence the signs of $w_{\pm 1}$ can be chosen such that
\begin{equation*}
w_{\pm 1}=\sqrt{\frac{s}{2}}\,\frac{1}{ps}\,(l\mp n),
\end{equation*}
or, from equation (\ref{eqn:e4}),
\begin{equation*}
w_{\pm 1}=\sqrt{\frac{s}{2}}\,(1-{e_4}^2)\,\frac{1}{A}(l\mp n).
\end{equation*}
Together with equation (\ref{eqn:omegasum}) and the equation for $dt$
obtained from (\ref{eqn:t}), this yields
\begin{align*}
i(\Omega_3+\omega_3)dt &= \frac{u-{e_4}^2}{1-{e_4}^2}\,w_{-1}\,\frac{du}{w}
+ i(u+1)\left(\frac{1}{C}-\frac{1}{A}\right)n\,dt\\
-i(\Omega_3-\omega_3)dt &= \frac{u-{e_4}^2}{1-{e_4}^2}\,w_{+1}\,\frac{du}{w}
+ i(u-1)\left(\frac{1}{C}-\frac{1}{A}\right)n\,dt.
\end{align*}
Substituting these expressions into (\ref{eqn:du-plus-iOmegas-dt}),
one obtains 
\begin{align*}
\frac{d\alpha}{\alpha} &= \frac{1}{2(u+1)}\left(w+
   \frac{u^2-{e_4}^2}{1-{e_4}^2}w_{-1}\right)\frac{du}{w}
+ \frac{i}{2}\left(\frac{1}{C}-\frac{1}{A}\right)n\,dt\\
\frac{d\beta}{\beta} &= \frac{1}{2(u-1)}\left(w-
   \frac{u^2-{e_4}^2}{1-{e_4}^2}w_{+1}\right)\frac{du}{w}
- \frac{i}{2}\left(\frac{1}{C}-\frac{1}{A}\right)n\,dt\\
\frac{d\gamma}{\gamma} &= \frac{1}{2(u-1)}\left(w+
   \frac{u^2-{e_4}^2}{1-{e_4}^2}w_{+1}\right)\frac{du}{w}
+ \frac{i}{2}\left(\frac{1}{C}-\frac{1}{A}\right)n\,dt\\
\frac{d\delta}{\delta} &= \frac{1}{2(u+1)}\left(w-
   \frac{u^2-{e_4}^2}{1-{e_4}^2}w_{-1}\right)\frac{du}{w}
- \frac{i}{2}\left(\frac{1}{C}-\frac{1}{A}\right)n\,dt.
\end{align*}
This implies equations (\ref{eqn:abgdints}) and the reduction to the
spherical case. 

Now assume that $A=C$. The assertions regarding the poles of the
logarithmic differentials away from $\infty$ follow
straightforwardly. Regarding the asymptotics at $\infty$, note that
near $\infty$ the logarithmic differentials are 
\begin{equation*}
\frac{du}{2(u\pm 1)}+\text{\it a holomorphic part}.
\end{equation*}
Introduce the parameter $\xi^2=1/(u\pm 1)$. It is well defined (up to
sign) and
holomorphic around $\infty$. One finds
that the singular part of the logarithmic differentials is
$d\xi/\xi$.
\end{proof}


\section{The trajectory of the tip of the top}\label{sec:tip}

In this section we examine the curve that the tip of the top
describes in the horizontal plane. This aspect of the top's motion is
particularly easy to observe experimentally. For example, in toy shops
one can buy tops that have a pencil lead for a tip or a felt tip
pen. Lord Kelvin did experiments with such pencil tipped tops, while
F. Klein used cogwheels from a clockwork which he spun on a sooty
glass plate. See \cite[pp. 619ff]{Klein:Buch} for pictures and a
discussion of the significant influence of friction and other
concomitants neglected in the mathematical model.

\subsection{Geometrical considerations}

Since $r$ is the unit vector in the direction of the top's axis, 
pointing away from the top's tip if attached to the center of gravity,
the curve traced by the tip of the top on the supporting plane is the
orthogonal projection of $-sr$ onto that plane. The following
proposition gives a formula for this curve in terms of the
Cayley-Klein parameters.
\begin{prop}
Identify the supporting plane with the complex number plane. Then the
curve traced by the tip of the top is
\begin{equation*}
c=2s\alpha\beta.
\end{equation*}  
\end{prop}
\begin{proof}
First, project the vector $-r$ stereographically from the north pole
of the unit sphere into the complex plane. 
The result is the image of $0$ under the M\"obius transformation
(\ref{eqn:Moeb}), i.e.\ $z=\beta/\delta$. Since $\bar{z}=-\gamma/\alpha$, 
the fromula for $c$ follows from equations
(\ref{eqn:inerse-stereo-proj}):
\begin{equation*}
c=s(\xi+i\eta)=s\,\frac{2z}{|z|^2+1}\,.
\end{equation*}
This proves the proposition.
\end{proof}

For the qualitative considerations below, it is useful to represent
the curve $c$ in
polar coordinates: $c=\rho e^{i\phi}$. Clearly,
$\rho=s\sqrt{1-u^2}$. The angle $\phi$ is not uniquely defined if
$\rho$ 
vanishes, i.e.\ if $u=\pm 1$. Otherwise, the following proposition
gives a formula for $d\phi$.
\begin{prop}
If $-1<u<1$, then 
\begin{equation}\label{eqn:dphi}
i\,d\phi=\frac{1}{2}\,\frac{u^2-{e_4}^2}{1-{e_4}^2}
\left(\frac{w_{-1}}{u+1}-\frac{w_{+1}}{u-1}\right)\frac{du}{w}\,.
\end{equation}
\end{prop}
\begin{proof}
The angle $\phi$ is the argument of the curve $z=2s\alpha\beta$. It follows
that $d\phi$ is the imaginary part of the logarithmic differential
$dc/c=d\alpha/\alpha+d\beta/\beta$. Formula (\ref{eqn:dphi}) then follows from
equations (\ref{eqn:abgdints}).
\end{proof}

\subsection{Loop, cusp, or wobbly arc?}

\begin{figure}
\includegraphics[width=0.3\textwidth]{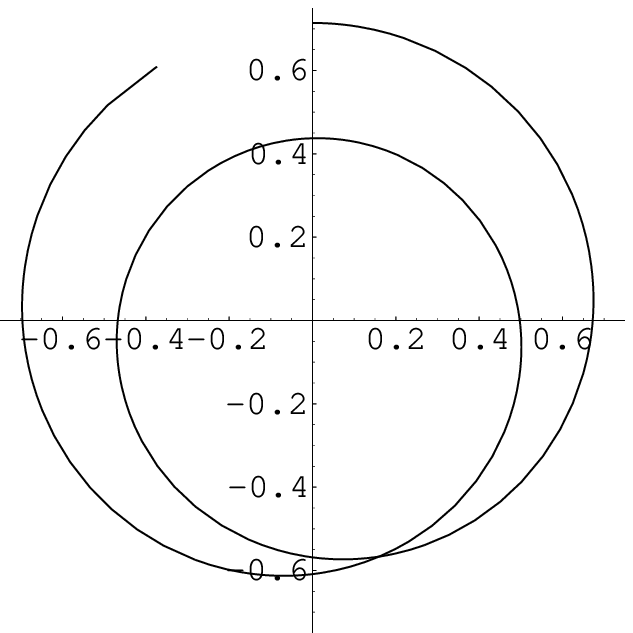}
\hfill
\includegraphics[width=0.3\textwidth]{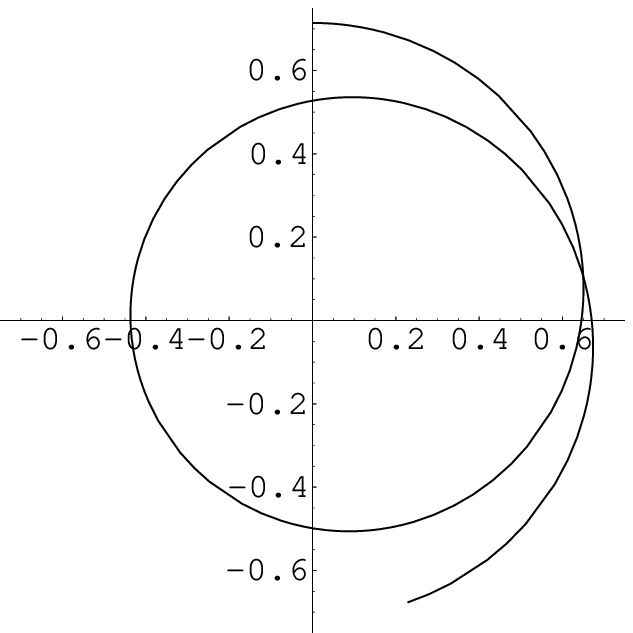}
\hfill
\includegraphics[width=0.3\textwidth]{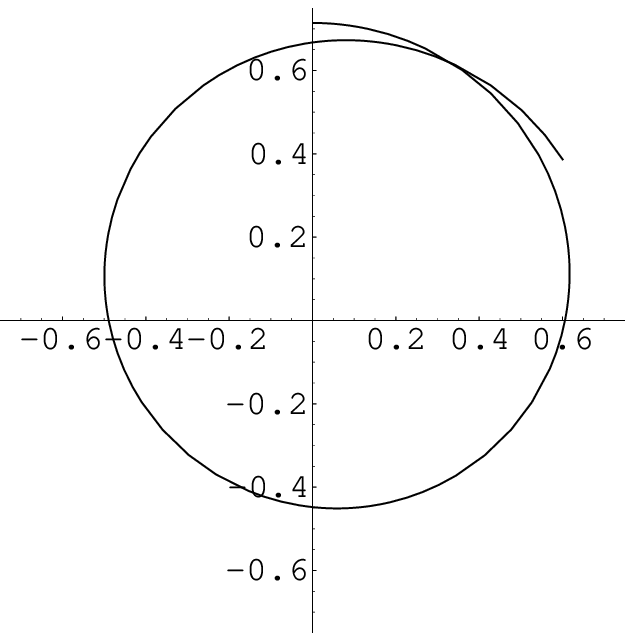}
\caption{The curve traced by the tip of the top. Here, $w_{-1}$ and
$w_{+1}$ have the same
sign, $s=1,\ e_1=0.7,\ e_2=0.9,\ {e_4}^2=2$ and  $e_3=1,\ 2,\ 100$
(left to right).}
\label{fig:tiptop-samesign}
\end{figure}
\begin{figure}
\includegraphics[width=0.3\textwidth]{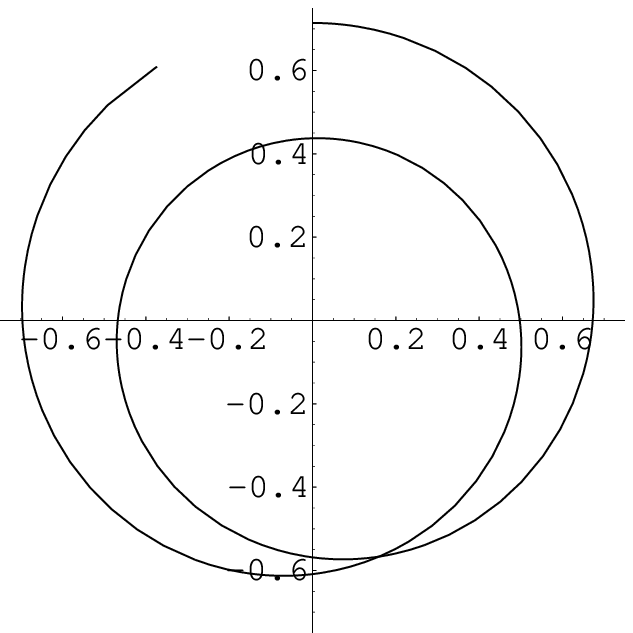}
\hfill
\includegraphics[width=0.3\textwidth]{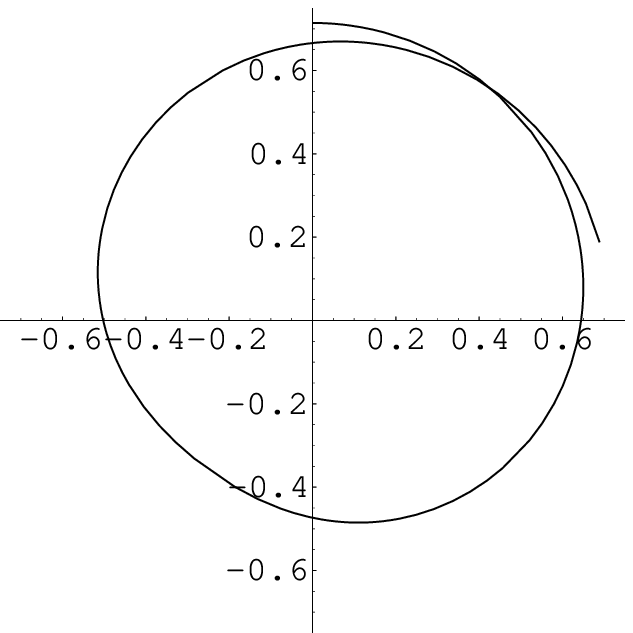}
\hfill
\includegraphics[width=0.3\textwidth]{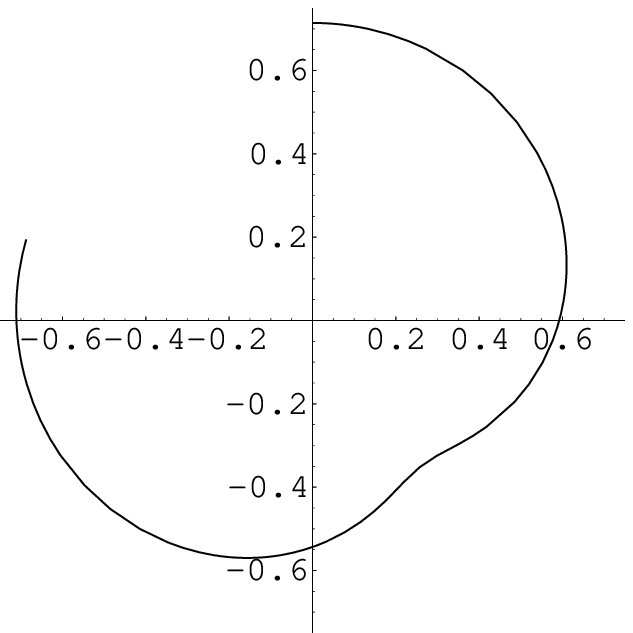} \\
\includegraphics[width=0.3\textwidth]{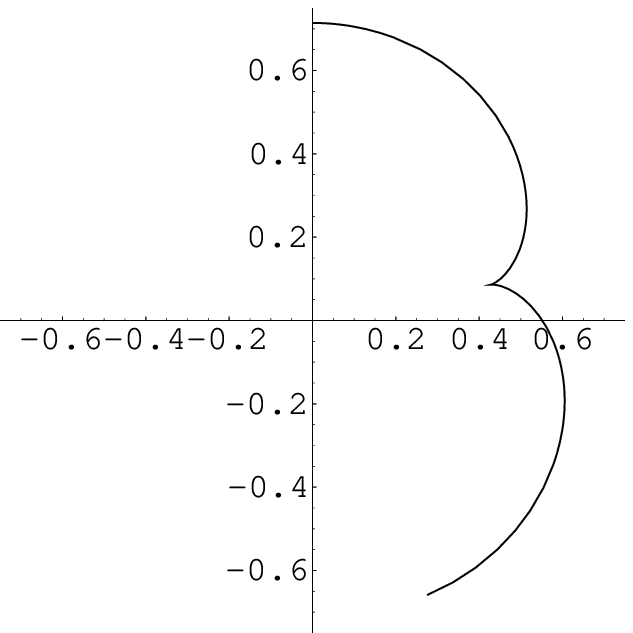}
\hfill
\includegraphics[width=0.3\textwidth]{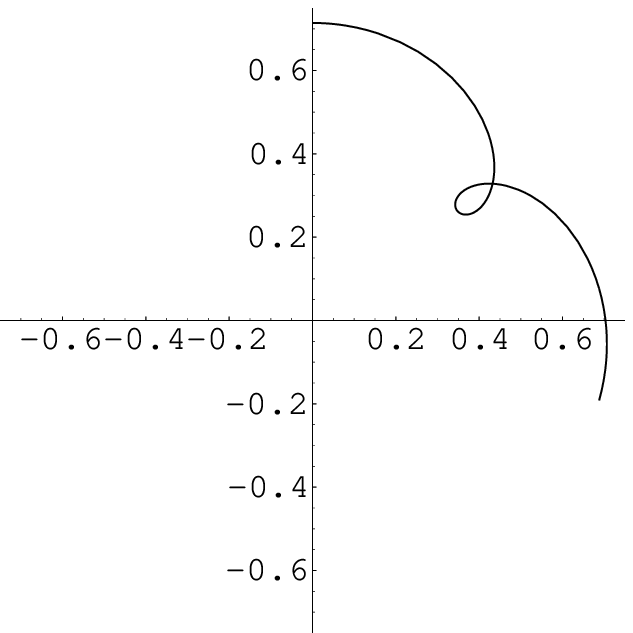}
\hfill
\includegraphics[width=0.3\textwidth]{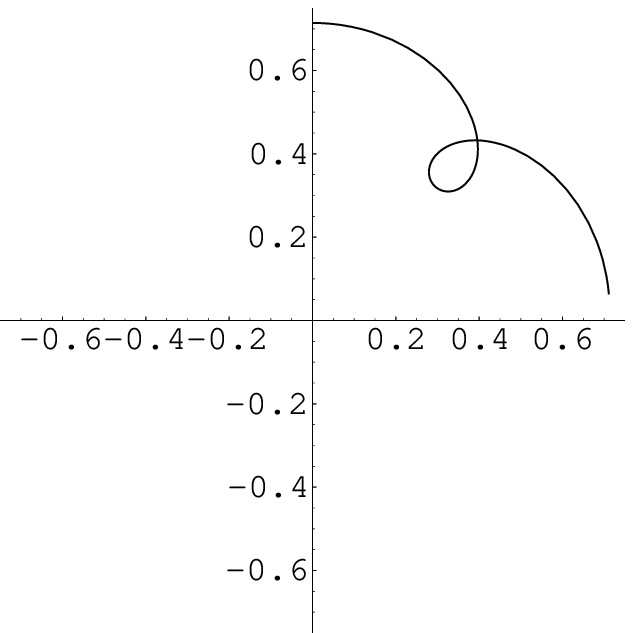}
\caption{The curve traced by the tip of the top. Here, $w_{-1}$ and $w_{+1}$ have different
signs, $s=1,\ e_1=0.7,\ e_2=0.9,\ {e_4}^2=2$ and  $e_3=1,\ 1.1,\
1.3,\ 1.85,\ 2.5,\ 3$  (top left to bottom right).} 
\label{fig:tiptop-diffrentsign}
\end{figure}
Figures \ref{fig:tiptop-samesign} and \ref{fig:tiptop-diffrentsign}
show different trajectories of the top's tip. Three qualitatively
different cases are clearly discernible. In Figure
\ref{fig:tiptop-samesign} and the top row of Figure
\ref{fig:tiptop-diffrentsign}, the curve is a smooth line circling the
origin. The first picture on the bottom row of Figure
\ref{fig:tiptop-diffrentsign} shows a cusp, and the following two show
loops. The next proposition gives a condition in
terms of $e_1$, $e_2$ and $e_3$ for the three cases to occur. 
\begin{prop}
Assume that $-1<e_1$ and $e_2<1$. If $w_{+1}$ and $w_{-1}$ have the
same sign, then the trajectory of the top's tip has neither loops nor
cusps. If $w_{+1}$ and $w_{-1}$ have different signs, then the curve
has loops if and only if
\begin{equation*}
\displaystyle{\frac{1-e_1 e_2}{e_2-e_1}} < e_3.
\end{equation*}
There are cusps if both sides are equal. 
\end{prop}
\begin{proof}
Loops occur, if $d\phi/du$ changes sign as $u$ moves from $e_1$ to
$e_2$. It follows from equation (\ref{eqn:dphi}) that this happens if 
$w_{-1}/(u+1)-w_{+1}/(u-1)=0$ has a solution for $u$ in the interval
$(e_1, e_2)$. That equation is equivalent to 
\begin{equation}\label{eqn:loopcondition}
\frac{u+1}{u-1}=\frac{w_{-1}}{w_{+1}}\,.
\end{equation}
Since the left hand side of this equation is smaller than zero for
$u\in (-1,1)$, this equation cannot be fulfilled if $w_{+1}$ and
$w_{-1}$ have the same sign. This proves the first part of the
proposition. 

Now assume that $w_{+1}$ and $w_{-1}$ have different signs. Then 
\begin{equation*}
\frac{w_{-1}}{w_{+1}}=-\sqrt{
\frac{(-1-e_1)(-1-e_2)(-1-e_3)}{(1-e_1)(1-e_2)(1-e_3)}
}.
\end{equation*}
In the interval $[e_1, e_2]$, the function $u\mapsto (u+1)/(u-1)$ is
continuously decreasing from $(e_1+1)/(e_1-1)$ to $(e_2+1)/(e_2-1)$.
Hence, equation (\ref{eqn:loopcondition}) leads to a $u\in (e_1, e_2)$
if 
\begin{equation*}
\frac{e_2+1}{e_2-1}<
-\sqrt{\frac{(-1-e_1)(-1-e_2)(-1-e_3)}{(1-e_1)(1-e_2)(1-e_3)}}<
\frac{e_1+1}{e_1-1}\,.
\end{equation*}
This is equivalent to 
\begin{equation*}
-\frac{1-e_3}{1+e_3}<\frac{1+e_1}{1-e_1}\,\frac{1-e_2}{1+e_2}<
-\frac{1+e_3}{1-e_3}.
\end{equation*}
The inequality on the right is always fulfilled, because the right
hand side is greater than one and the left hand side smaller. The
inequality on the left
is equivalent to $(1-e_1 e_2)/(e_2-e_1)<e_3$.
This proves the condition for loops. If $(1-e_1 e_2)/(e_2-e_1)=e_3$,
then the zero of $d\phi/du$ occurs at $u=e_2$, hence there are cusps.
\end{proof}


\section{The degenerate cases}

If two zeroes of the polynomial $R(u)$ coincide, the hyperelliptic
integrals of the Sections \ref{sec:t-int} and \ref{sec:abgd-ints}
degenerate to elliptic ones. In 
those cases, it is possible to solve the system in terms of elliptic
functions. This will be done in the following
sections. Since $-1\leq e_1\leq e_2\leq 1\leq e_3$ and $e_4>1$,
only the following cases have to be considered: $e_1=e_2$, $e_2=e_3$
and $e_3=e_4$. Furthermore, the limit case $e_4=-e_4=\infty$ is
examined, because it is of special interest: In this case the toy top
becomes the Lagrange top. 


\subsection{The case $e_3=e_4$}\label{sec:e3=e4}

Suppose that $e_3=e_4$. Then the hyperelliptic curve $C$ degenerates
to the elliptic curve $C_{e_3=e_4}$, given by $w^2=R(u)$, where
\[ R(u)=4(u-e_1)(u-e_2)(u+e_3)\, .
\]
The leading coefficient of $R$ is chosen to be 4 in accordance with
the Weierstrass normalization of elliptic curves. This will simplify
the integration of the elliptic integrals in terms of Weierstass
elliptic functions $\wp, \zeta$ and $\sigma$. The branchpoints are
$e_1, e_2, -e_3$, and $\infty$.  

The $t$-integral and the integrals (\ref{eqn:abgdints}) become
\begin{equation}\label{eqn:degtint}
 t= -i\sqrt{2s}\int (u+e_3)\frac{du}{w}\,,
\end{equation}
and
\begin{equation}\label{eqn:zander}
\begin{split}
\log \alpha &= \int \frac{1}{2(u+1)}\left(w+
   \frac{u+e_3}{-1+e_3}\,w_{-1}\right)\frac{du}{w}\\
\log \beta &= \int \frac{1}{2(u-1)}\left(w-
   \frac{u+e_3}{1+e_3}\,w_{+1}\right)\frac{du}{w}\\
\log \gamma &= \int \frac{1}{2(u-1)}\left(w+
   \frac{u+e_3}{1+e_3}\,w_{+1}\right)\frac{du}{w}\\
\log \delta &= \int \frac{1}{2(u+1)}\left(w-
   \frac{u+e_3}{-1+e_3}\,w_{-1}\right)\frac{du}{w}\,.
\end{split}
\end{equation}
As before, $w_{\pm 1}$ denotes one of the two values of $w$ on the
curve $C_{e_3=e_4}$ over $u=\pm 1$, i.e
\begin{equation*}
w_{+1}=\pm\sqrt{R(1)}\quad\text{and}\quad w_{-1}=\pm\sqrt{R(-1)}.
\end{equation*}


\subsubsection{Uniformization of the elliptic curve}
\label{sec:uniformization}

The system will be solved by using an Abelian integral of the first
kind to uniformize the the elliptic curve $C_{e_3=e_4}$ and then
expressing $t$ and $\alpha, \beta, \gamma, \delta$ as functions of the
uniformizing variable. 

\begin{figure}
\begin{center}
\parbox{0.54\textwidth}{
   \includegraphics[width=0.5\textwidth]{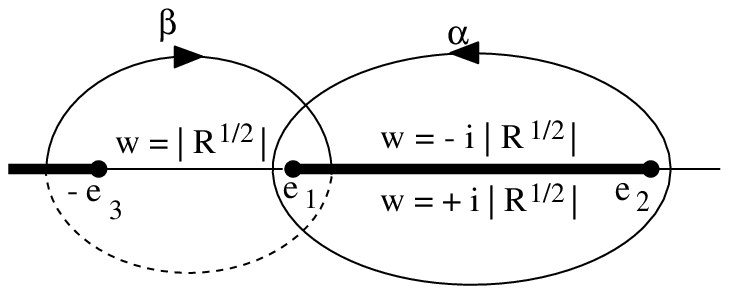}
}
\parbox{0.41\linewidth}{
   \includegraphics[width=0.41\textwidth]{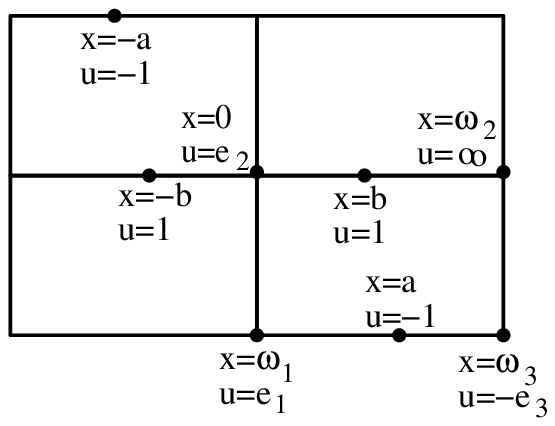}
}
\end{center}
\caption{Left: The top sheet of the elliptic curve $C_{e_3=e_4}$. 
Right: A fundamental domain in the $x$-plane.}
\label{fig:curve2}
\end{figure}
The left part of Figure \ref{fig:curve2} shows the top sheet of
the curve $C_{e_3=e_4}$ with cuts and a normal cycle basis
$\alpha,\beta.$ Let
$\omega_1$ and $\omega_2$ be the corresponding half periods:
\begin{equation*}
2\omega_1=\int_\alpha\frac{du}{w},\quad 
   2\omega_2=\int_\beta\frac{du}{w}\,.
\end{equation*}
Note that $\omega_1$ lies on the negative imaginary axis and that $\omega_2$
is positive. Let the uniformizing variable be given by 
\begin{equation}\label{eqn:x}
x=\int_{u=e_2}^{u=u(x)}\frac{du}{w}.
\end{equation}
One obtains $u$ and $w$
as doubly periodic functions of $x$ with periods $2\omega_1$ and
$2\omega_2$, namely
\begin{equation}\label{eqn:uofx}
u(x)=\wp(x-\omega_2)+\frac{1}{3}(e_1+e_2-e_3)
\end{equation}
and 
\begin{equation}\label{eqn:wofx}
w(x)=\wp'(x-\omega_2).
\end{equation}
The additive constant $(e_1+e_2-e_3)/3$ in (\ref{eqn:uofx}) comes
from the fact that the branchpoints are not centered around zero as
required by the Weierstrass normalization, and the $\wp$-function is
shifted by $\omega_1$ because the integral in (\ref{eqn:x}) does not
start in $\infty$ but $e_2$. 

The right part of Figure \ref{fig:curve2}
shows a fundamental rectangle in the $x$-plane. The points $x=\pm a$
and $x=\pm b$ will be of importance in Section
\ref{sec:abgdformulas} and are defined as follows: The point $a$ in
the $x$-plane is supposed to correspond to the
point $(u,w)=(-1,w_{-1})$ on $C_{e_3=e_4}$, and $b$ is supposed to
correspond to $(u,w)=(+1,w_{+1})$. 
For example, let
\begin{equation*}
a=\pm\left(\omega_1+\int_{-1}^{e_1}\sqrt{R(u)}\,du\right)
\end{equation*}
and
\begin{equation*}
b=\pm\int_{e_2}^{1}\sqrt{R(u)}\,du,
\end{equation*}
where the signs are chosen according to whether the points
$(u,w_{\pm1})$ lie in the bottom or top sheet. (The figure shows the
case where both $x=a$ and $x=b$ correspond to points in the bottom
sheet. Otherwise exchange $+a$ with $-a$, or $+b$ with $-b$,
respectively.) 

Regarding the toy top we adopt the convention that during its motion,
the corresponding point on $C_{e_3=e_4}$ moves on a path homotopic to
the cycle $\alpha$. (This choice determines the sign on the right hand
side of equation (\ref{eqn:degtint})). The corresponding point
in the $x$-plane then moves 
on the imaginary axis in the negative direction.


\subsubsection{Solution for time as function of the
uniformizing variable}

We will now pull back the $t$-integral (\ref{eqn:degtint}) to the
$x$ plane and solve it. The following proposition gives the resulting
formula for $t$ as function of $x$. The initial condition which is
chosen means that $u=e_2$ at $t=0$, i.e.\ the axis of the top is
initially in its most upright position. 
\begin{prop}
Assume that $t=0$ for $x=0$. Then
equation (\ref{eqn:degtint}) implies
\begin{equation*}
 t=i\sqrt{2s}\big(\zeta(x-\omega_2)+\zeta(\omega_2)
-\frac{1}{3}(e_1+e_2+2e_3)\,x\big).
\end{equation*}
\end{prop}
\begin{proof}
Equations (\ref{eqn:uofx}) and (\ref{eqn:wofx}) imply that the
pullback of the holomorphic differential $du/w$ to the $x$-plane is
$dx$. The $t$-integral therefore becomes
\begin{equation*}
 t=-i\sqrt{2s}\int\big(\wp(x-\omega_2)+\frac{1}{3}(e_1+e_2+2e_3)\big)dx.
\end{equation*}
Since $\zeta'=-\wp$, this implies
\begin{equation*}
 t=i\sqrt{2s}\big(\zeta(x-\omega_2)
-\frac{1}{3}(e_1+e_2+2e_3)x+\mbox{\em const.}\big).
\end{equation*}
The constant has to be $\zeta(\omega_2)$ for $t$ to vanish at
$x=0$.
\end{proof}


\subsubsection{Solution for the Cayley-Klein parameters as functions
of the uniformizing variable}
\label{sec:abgdformulas}

The following proposition gives formulas for $\alpha, \beta, \gamma, \delta$ as
functions of $x$. But first, a few words have to be said about the
initial conditions $\alpha_0, \beta_0, \gamma_0, \delta_0$ for $t=x=0$. 

It is no essential restriction to assume that $\alpha_0$ is real and
nonnegative and therefore equals $\delta_0$, and that $\gamma_0$ is purely
imaginary with nonnegative imaginary part, such that
$\beta_0=\gamma_0$. For one may always achieve this by 
suitable rotations of the system about the $z$-Axis and the top's
symmetry axis.

\begin{prop}
Under the above initial conditions, one obtains 
\begin{equation}\label{eqn:degabgd}
\begin{split}
\alpha=k_1 e^{l_1 x}\frac{\sigma(x-a)}{\sigma(x-\omega_2)}, \\
\beta=k_2 e^{l_2 x}\frac{\sigma(x+b)}{\sigma(x+\omega_2)}, \\
\gamma=k_3 e^{l_3 x}\frac{\sigma(x-b)}{\sigma(x-\omega_2)}, \\
\delta=k_4 e^{l_4 x}\frac{\sigma(x+a)}{\sigma(x+\omega_2)},
\end{split}
\end{equation}
where
\begin{equation}\label{eqn:ks}
\begin{split}
k_1 &= k_4 =
        \sqrt{\frac{1+e_2}{2}}\;\frac{\sigma(\omega_2)}{\sigma(a)},\\
k_2 &= k_3 =
   \,i\,\sqrt{\frac{1-e_2}{2}}\;\frac{\sigma(\omega_2)}{\sigma(b)},
\end{split}
\end{equation}
and
\begin{equation}\label{eqn:ls}
\begin{split}
l_1 &= -l_4 = \frac{w_{-1}}{2(-1+e_3)}+\zeta(a-\omega_2),\\
l_3 &= -l_2 = \frac{w_{+1}}{2(1+e_3)}+\zeta(b-\omega_2).
\end{split}
\end{equation}
\end{prop}
\begin{proof}
Consider the asymptotic behavior of the logarithmic differential
of $\alpha$:
\[ \frac{d\alpha}{\alpha}=\frac{1}{2(u+1)}\left(w+
   \frac{u+e_3}{-1+e_3}\,w_{-1}\right)\frac{du}{w}.
\]
It has only two simple poles on $C_{e_3=e_4}$, one at
$(u,w)=(-1,w_{-1})$ and the other at $\infty$. For
$(u,w)\rightarrow(-1,w_{-1})$, the asymptotic expansion is
\[\frac{d\alpha}{\alpha}=\left(\frac{1}{u+1}+O(1)\right)du.
\]
To calculate the asymptotic expansion around the branchpoint
$\infty$, introduce the local parameter $\xi^2=1/u$. The result is 
\[\frac{d\alpha}{\alpha}= \left(-\frac{1}{\xi}+O(1)\right)d\xi.
\]
Since the residues of the logarithmic differential of $\alpha$ are $1$
and $-1$, respectively, $\alpha$ itself is not branched on $C_{e_3=e_4}$,
but has a zero at $(u,w)=(-1,w_{-1})$ and a simple pole at
$\infty$. It is a multiply valued function, because the additive
periods of the 
logarithmic differential lead to multiplicative periods of $\alpha$. It
follows that, in terms of the uniformizing variable $x$, $\alpha$ is of
the form 
\begin{equation*}
\alpha=k_1 e^{l_1 x}\frac{\sigma(x-a)}{\sigma(x-\omega_2)}.
\end{equation*}
Similarly, one obtains the other equations (\ref{eqn:degabgd}).

The constants $k_1,\ldots,k_4$ and $l_1,\ldots,l_4$ are determined by
the initial conditions. Note first that $\alpha_0=\delta_0$
implies $k_1=k_4$ and $\beta_0=\gamma_0$ implies $k_2=k_3$. Just substitute
$x=0$ into equations (\ref{eqn:degabgd}) and observe that $\sigma$ is
an odd function. 

Further, since $u=e_2$ for $x=0$, it follows from equations
(\ref{eqn:up1um1}) that
\begin{equation*}
e_2+1=2\alpha_0\delta_0=2{k_1}^2\frac{\sigma^2(a)}{\sigma^2(\omega_2)}
\end{equation*}
and
\begin{equation*}
e_2-1=2\beta_0\gamma_0=2{k_2}^2\frac{\sigma^2(b)}{\sigma^2(\omega_2)}.
\end{equation*}
From this one obtains equations (\ref{eqn:ks}). The signs of $k_1$ and
$k_2$ are determined so that $\alpha_0$ is positive
and $\gamma_0$ has a positive imaginary part.

Concerning the constants $l_1,\ldots,l_4$, note first that
   \begin{equation*}
   \frac{u+1}{2}=\alpha\delta=k_1 k_4e^{(l_1+l_4)x}
   \frac{\sigma(x-a)\sigma(x+a)}{\sigma(x-\omega_2)\sigma(x+\omega_2)}
   \end{equation*}
is a doubly periodic function of $x$. This implies $l_4=-l_1$ since
the quotient of $\sigma$-functions is already doubly
periodic. Equally, considering $\beta\gamma=(u-1)/2$ yields $l_2=-l_3$. 

From the logarithmic derivative of $\alpha$ with respect to $x$
\begin{equation*}
\frac{d\log\alpha}{dx}=
   l_1+\frac{\sigma'(x-a)}{\sigma(x-a)}
      -\frac{\sigma'(x-\omega_2)}{\sigma(x-\omega_2)}
\end{equation*}
one obtains
\begin{equation*}
l_1=\frac{d\log\alpha}{dx}-\zeta(x-a)+\zeta(x-\omega_2),
\end{equation*}
since $\sigma'/\sigma=\zeta$. The first of equations then
(\ref{eqn:ls}) follows from 
\begin{equation*}
\lim_{x\rightarrow a}\left(\frac{d\log\alpha}{dx}-\zeta(x-a)\right) =
   \frac{w_{-1}}{2(-1+e_3)}\,.
\end{equation*}
Since $\zeta(x-a)=1/(x-a)+O(x-a)$ as $x\rightarrow a$, one has to show
that 
\begin{equation}\label{eqn:alphaasymp}
\frac{d\log\alpha}{dx}=\frac{1}{x-a}+\frac{w_{-1}}{2(-1+e_3)}+O(x-a).
\end{equation}
From the first equation (\ref{eqn:zander}) and $dx=du/w$
one obtains 
\begin{equation*}
\frac{d\log\alpha}{dx}=
   \frac{1}{2(u+1)}
     \left(w+\frac{u+e_3}{-1+e_3}\,w_{-1}\right).
\end{equation*}
Since $w=du/dx$, this yields the expansion
\begin{align}\label{eqn:flunder}
\frac{d\log\alpha}{dx} 
 &= \frac{1}{2(u+1)}\left(
   \left.\frac{du}{dx}
      \right|_{\makebox[0pt][l]{$\scriptstyle x=a$}}+
   (x-a)\left.\frac{d^2u}{dx^2}
      \right|_{\makebox[0pt][l]{$\scriptstyle x=a$}}+
   O(x-a)^2 +
   \left(1+\frac{u+1}{-1+e_3}\right)w_{-1}\right) \notag\\
 &= \frac{1}{u+1}\left.
   \frac{du}{dx}\right|_{x=a}+
   \frac{1}{2}\left.\frac{\,\frac{d^2u}{dx^2}\,}
                   {\frac{du}{dx}}\right|_{x=a}
    +
   \frac{w_{-1}}{2(-1+e_3)}+O(x-a)\,.
\end{align}
The last equality uses 
\begin{equation}\label{eqn:diffquot}
\frac{x-a}{u+1}=\frac{1}
   {\left.\frac{du}{dx}\right|_{\makebox[0pt][l]{$\scriptstyle x=a$}}}
+O(x-a).
\end{equation}
The asymptotical expansion (\ref{eqn:alphaasymp}) follows from
(\ref{eqn:flunder}) and
\begin{equation*}
\frac{1}{x-a}=
\frac{1}{u+1}\left.
   \frac{du}{dx}\right|_{x=a}+
   \frac{1}{2}\left.\frac{\,\frac{d^2u}{dx^2}\,}
                   {\frac{du}{dx}}\right|_{x=a}+O(x-a).
\end{equation*}
To see this last equation, just take the Taylor expansion
\begin{equation*}
u+1=\left.\frac{du}{dx}\right|_{x=a}(x-a)+
    \frac{1}{2}\left.\frac{d^2u}{dx^2}\right|_{x=a}(x-a)^2+
    O(x-a)^3,
\end{equation*}
divide by $(x-a)$ and $(u+1)$, and use (\ref{eqn:diffquot}) again.
The equation for $l_3$ is derived analogously.
\end{proof}


\subsection{The case of regular precession}\label{sec:reg-prec}

If a spinning top moves in such a way that the angle between the top's
axis and the vertical is constant, its motion is called regular
precession. It then follows from the symmetries of the system that
the top spins around its axis with constant speed,
while the axis precedes uniformly around the vertical axis. Notably,
the solutions for $t$ and the Cayley-Klein parameters in terms of
hyperelliptic integrals fail in the case of regular precession because
the path of integration degenerates to a point. One would have to look
at the system before it is completely reduced to analyze this case. It
turns out that regular precession occurs if $R(u)$ has a
double zero at the constant value for $u$. There are two qualitatively
different ways this can happen. If $e_1=e_2$, small
perturbations of the system will lead to small nutation. The regular
precession is stable. But if
$e_2=e_3=1$, small perturbations will lead to nutation between
$u=e_1$ and $u=e_2$. The regular precession is unstable. See
\cite[pp. 278ff]{Klein:Buch} for an analogous
discussion of the stability of regular precession in the case of
Lagrange's top.


\subsection{The aperiodic case: $e_2=e_3=1$}\label{sec:aper}

If $e_2=e_3$, there are two possible types of motion. 
Either $u$ is constantly 1. This is the unstable case of regular
precession mentioned above. Or $u$ moves from 1
to $e_1$ and back. This case will be considered in this
section. We proceed similarly as in Section \ref{sec:e3=e4}.

The hyperelliptic curve $C$ degenerates to the elliptic
curve $C_{\text{aper}}$ given by $w^2=R(u)$ with
\begin{equation*}
R(u)=4(u-e_1)(u^2-{e_4}^2).
\end{equation*}
Branchpoints are $e_1$, $\pm e_4$ and $\infty$.
The $t$-integral becomes
\begin{equation}\label{eqn:apertint}
t = i\,\sqrt{2s} \int \frac{(u^2-{e_4}^2)}{u-1}\,\frac{du}{w}.
\end{equation}
It is singular at $u=1$. That is why this is the aperiodic case.
The integrals (\ref{eqn:abgdints}) become
\begin{equation}\label{eqn:aperabgdints}
\begin{split}
\log\alpha&= \int\frac{1}{2(u+1)}
   \left( w-\frac{2}{u-1}\,\frac{u^2-{e_4}^2}{1-{e_4}^2}w_{-1}
   \right) \frac{du}{w} \\
\log\beta &= \int\frac{du}{2(u-1)} \\
\log\gamma &= \int\frac{du}{2(u-1)} \\
\log\delta &= \int\frac{1}{2(u+1)}
   \left( w+\frac{2}{u-1}\,\frac{u^2-{e_4}^2}{1-{e_4}^2}w_{-1}
   \right) \frac{du}{w}\;.
\end{split}
\end{equation}


\subsubsection{Uniformization of the elliptic curve}

\begin{figure}
\begin{center}
\parbox{.5\textwidth}{
   \includegraphics[width=0.5\textwidth]{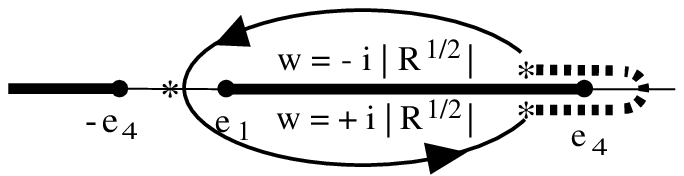}
}
\parbox{.4\textwidth}{
   \includegraphics[width=0.4\textwidth]{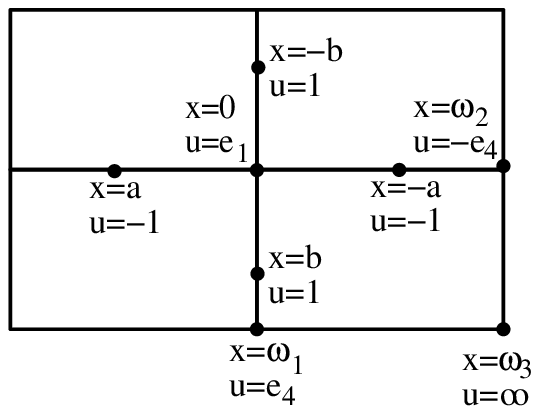}
}
\caption{Left: The top sheet of the elliptic curve
$C_{\text{aper}}$. Right: A fundamental domain in the $x$-plane.}
\label{fig:curve3}
\end{center}
\end{figure}

Figure \ref{fig:curve3} shows one sheet of $C_{\text{aper}}$, which
will be called the top sheet.
Proceed as in Section \ref{sec:uniformization}, but let the curves $c$
start in $(u,w)=(e_1,0)$. For a normalized cycle basis, take a path
$\alpha$ in the top sheet which encircles $e_1$ and $e_4$ in the
counterclockwise direction, and choose $\beta$ accordingly. Let $\omega_1$
and $\omega_2$ be the corresponding half periods and let
$\omega_3=\omega_1+\omega_2$.
One obtains the uniformization
\begin{equation*}
u=\wp(x-\omega_3)+\frac{e_1}{3}
\end{equation*}
and $w=\wp'(x-\omega_3)$.

The right part of Figure \ref{fig:curve3} shows a
fundamental domain in the $x$-plane. Again, the points above $u=\pm 1$
will be important. In the figure, they are marked by asterisks. Let
$a$ and $b$ be points in the $x$-plane corresponding to
$(u,w)=(-1,w_{-1})$ and $(u,w)=(+1,w_{+1})$, for example, 
\begin{equation*}
a=\pm\int_{-1}^{e_1}\sqrt{R(u)}\,du
\quad
\text{and}
\quad
b=\pm i\int_{e_1}^{1}\sqrt{-R(u)}\,du.
\end{equation*}
(In the figure they are drawn for the case in which $(-1,w_{-1})$ lies
in 
the top sheet and $w_{+1}$ has positive imaginary part.)

Let the path of integration corresponding to the motion of the top be
homotopic the the path with arrows drawn in the figure. 


\subsubsection{Solution for time as function of the uniformizing
variable}

The integrand of the $t$-integral has simple poles over $u=1$ which
leads to logarithmic type singularities for $t$. The Riemann surface
of the function $t$ will therefore be an infinitely sheeted cover of
$C_{\text{aper}}$ with the two branchpoints $(u,w)=(1, \pm w_{+1})$. 
Place a cut on $C_{\text{aper}}$ along the dotted line in
Figure \ref{fig:curve3} and cut the $x$-plane along corresponding
lines. Then our path of integration on $C_{\text{aper}}$ does not
cross the cut, so that we can consider only one branch of the function
$t$ on the cut elliptic curve, or the cut $x$-plane, respectively.
\begin{prop}
Suppose that $t=0$ at $x=0$. Then
\begin{align*}
\frac{t}{i\sqrt{2s}}=
   &-\zeta(x-\omega_3)-\zeta(\omega_3)+\left(1+\frac{e_1}{3}\right)x\\
   &+\frac{1-{e_4}^2}{w_{+1}}\left(
      \log\left(-\frac{\sigma(x-b)}{\sigma(x+b)}\right)
      +2\big(\zeta(\omega_3)-\zeta(b-\omega_3)\big) x
   \right),
\end{align*}
where that branch of the function 
\begin{equation*}
\log\left(-\frac{\sigma(x-b)}{\sigma(x+b)}\right)
\end{equation*}
is chosen that is zero for $x=0$.
\end{prop}
\begin{proof}
Put equation (\ref{eqn:apertint}) in the form
\begin{equation}\label{eqn:tpolydiv}
t=i\sqrt{2s}\int\left(u+1+\frac{1-{e_4}^2}{u-1}\right)\frac{du}{w} 
\end{equation}
and consider separately the integrals 
\begin{equation*}
I_1=\int\frac{1}{u-1}\;\frac{du}{w}
\end{equation*}
and
\begin{equation*}
I_2=\int (u+1)\,\frac{du}{w}\,.
\end{equation*}
Substituting the uniformizing variable in the first integral, we get
\begin{align*}
I_1 &=\int\frac{dx}{\wp(x-\omega_3)-\wp(b-\omega_3)} \\
    &=\frac{1}{\wp'(b-\omega_3)}
         \int\frac{\wp'(b-\omega_3)\,dx}{\wp(x-\omega_3)-\wp(b-\omega_3)}. 
\end{align*}
Using the formula
\begin{equation*}
\frac{\wp'(\eta)}{\wp(\xi)-\wp(\eta)}=
\zeta(\xi-\eta)-\zeta(\xi+\eta)+2\zeta(\eta),
\end{equation*}
one obtains
\begin{align*}
I_1 
&= \frac{1}{\wp'(b-\omega_3)}
    \int\big(\zeta(x-b)-\zeta(x+b-2\omega_3)+2\zeta(b-\omega_3)\big) dx \\ 
&=  \frac{1}{\wp'(b-\omega_3)}
      \int\big(\zeta(x-b)-\zeta(x+b)+
      2\zeta(\omega_3)+2\zeta(b-\omega_3)\big) dx.
\end{align*}
The last equality follows from
$\zeta(\xi+2\omega_k)=\zeta(\xi)+2\zeta(\omega_k).$
Now the last integral can easily be solved since
$\zeta(\xi)={\sigma'(\xi)}/{\sigma(\xi)}$.
Note that
\begin{equation*}
\wp'(b-\omega_3)=w_{+1}
\end{equation*}
to obtain
\begin{equation*}
I_1=\frac{1}{w_{+1}}\left(
   \log\left(\frac{\sigma(x-b)}{\sigma(x+b)}\right)
   +2\big(\zeta(\omega_3)+\zeta(b-\omega_3)\big)x+\mbox{\em const.}
\right).
\end{equation*}
Requiring that $I_1=0$ for $x=0$, one obtains
\begin{equation}\label{eqn:I1}
I_1=\frac{1}{w_{+1}}\left(
   \log\left(-\frac{\sigma(x-b)}{\sigma(x+b)}\right)
   +2\big(\zeta(\omega_3)+\zeta(b-\omega_3)\big)x \right),
\end{equation}
where that branch of the logarithmic term is chosen that vanishes for
$x=0$. 
The other integral is easier to deal with:
\begin{align*}
I_2 &= \int\big(\wp(x-\omega_3)+1+e_1/3\big) \\
    &= -\zeta(x-\omega_3)+(1+e_1/3)x+\mbox{\em const.},
\end{align*}
since $\wp=-\zeta'$. If the constant is chosen to make $I_2$
vanish for $x=0$, this becomes
\begin{equation}\label{eqn:I2}
I_2 = -\zeta(x-\omega_3)-\zeta(\omega_3)+\left(1+\frac{e_1}{3}\right)x.
\end{equation}
The formula for $t$ follows from equations
(\ref{eqn:tpolydiv}), (\ref{eqn:I1}) and (\ref{eqn:I2}).
\end{proof}


\subsubsection{Solution for the Cayley-Klein parameters}

The functions $\alpha,\beta,\gamma,\delta$ are also branched at the points above
$u=1$. So we apply the same cut as in the previous section and have to
choose branches. As in Section
\ref{sec:e3=e4} we will assume that the initial conditions $\alpha_0$ and
$\delta_0$ are real (and therefore equal) and positive, and that $\beta_0$
and $\gamma_0$ are imaginary (and therefore equal) with positive
imaginary part. 

\begin{prop}
If the initial conditions are chosen as explained above, then
\begin{equation*}
\begin{split}
\alpha &= k e^{lx-i\pi p} \frac{\sigma(x-a)}{\sigma(x-\omega_3)}
   \left(\frac{\sigma(x+b)}{\sigma(x-b)}\right)^p \\
\gamma &= \beta = i\sqrt{\frac{1-u}{2}}\\
\delta &= k e^{-lx+i\pi p} \frac{\sigma(x+a)}{\sigma(x+\omega_3)}
   \left(\frac{\sigma(x-b)}{\sigma(x+b)}\right)^p,
\end{split}
\end{equation*}
where
\begin{gather*}
p=-\frac{1}{2}\,\frac{w_{-1}}{w_{+1}},\\
k=\sqrt{\frac{1+e_1}{2}}\,\frac{\sigma(\omega_3)}{\sigma(a)},\\
l=\frac{w_{-1}}{4}-\frac{w_{-1}}{1-{e_4}^2}+\zeta(a-\omega_3)
  -p\big(\zeta(a+b)+\zeta(a-b)),
\end{gather*}
and those branches of the multiply valued functions \mbox{$(\sigma(x+b)/\sigma(x-b))^p$} and \\
\mbox{$(\sigma(x-b)/\sigma(x+b))^p$} 
are chosen which take the
values $e^{i\pi p}$ and $e^{-i\pi p}$ at $x=0$.
\end{prop}
\begin{proof}
The formula for $\beta$ and $\gamma$ follows elementarily from the
corresponding integrals (\ref{eqn:aperabgdints}) and the initial
conditions. 
Now consider the logarithmic differential of $\alpha$ from equations
(\ref{eqn:aperabgdints}):
\begin{equation*}
\frac{d\alpha}{\alpha}=
\left( \frac{1}{2(u+1)}-\frac{1}{(u-1)(u+1)}\,\frac{w_{-1}}{w}
   -\frac{1}{1-{e_4}^2}\,\frac{w_{-1}}{w}
\right)du .
\end{equation*}
As one can read off from this expression, the differential has three
simple poles away from infinity, one at
$(u,w)=(-1,w_{-1})$ with residue $1$ and two more at
$(u,w)=(+1,\pm w_{+1})$ with residues $\mp w_{-1}/2w_{+1}$. 
Furthermore, the term $du/2(u+1)$ contributes one
more simple pole at infinity with residue $-1$. 
\begin{table}
\caption{The simple poles of $d\alpha/\alpha$ and $d\delta/\delta$ with
residues.}
\label{tab:aperpoles}
\small
\begin{center}
\begin{tabular}{l|c|c|c|c}
\multicolumn{5}{c}{$d\alpha/\alpha$ } \\
$(u,w)=$&
$(-1,w_{-1})$&$\quad\infty\quad$& $(+1,+w_{+1})$ & $(+1,-w_{+1})$\\
\hline
Residue: &
1            & -1              &$-w_{-1}/2w_{+1}$&$w_{-1}/2w_{+1}$\\ 
\\[-1ex]
\multicolumn{5}{c}{$d\delta/\delta$ } \\
$(u,w)=$&
$(-1,-w_{-1})$&$\quad\infty\quad$& $(+1,+w_{+1})$ & $(+1,-w_{+1})$\\
\hline
Residue: &
1            & -1              &$w_{-1}/2w_{+1}$&$-w_{-1}/2w_{+1}$
\end{tabular}
\end{center}
\end{table}
The location of the poles with their residues is summarized in Table
\ref{tab:aperpoles}. It also lists the poles of the logarithmic
differential of $\beta$, which are obtained similarly. Now the poles
with residue $\pm 1$ lead to zeroes and poles of $\alpha$ and $\delta$,
while the other poles give rise to branchpoints. It follows that, as
functions of $x$, $\alpha$ and $\delta$ are of the form
\begin{equation}\label{eqn:aperform}
\begin{split}
\alpha &= k_1 e^{l_1 x} \frac{\sigma(x-a)}{\sigma(x-\omega_3)}
   \left( \frac{\sigma(x+b)}{\sigma(x-b)} \right)^p  \\
\beta &= k_2 e^{l_2 x} \frac{\sigma(x+a)}{\sigma(x+\omega_3)}
   \left( \frac{\sigma(x-b)}{\sigma(x+b)} \right)^p. 
\end{split}
\end{equation}

Choose the branches of these multiply valued functions as explained in
the proposition. 
Then the values
\begin{align*}
k_1 &= \sqrt{\frac{1+e_1}{2}}\,\frac{\sigma(\omega_3)}{\sigma(a)}\,
       e^{-i\pi p} \\
k_2 &= \sqrt{\frac{1+e_1}{2}}\,\frac{\sigma(\omega_3)}{\sigma(a)}\,
       e^{+i\pi p}
\end{align*}
follow from the initial condition $\alpha_0=\delta_0=\sqrt{(1+e_1)/2}$. 

Since $2\alpha\delta=u+1$ is a doubly periodic function of $x$, it follows
that $l_2=-l_1$. It is left to show that $l_1=l$ as given in the
proposition. Since from equation (\ref{eqn:aperform})
\begin{equation*}
\frac{d\log\alpha}{dx}=l_1+\zeta(x-a)-\zeta(x-\omega_3)
   +p\big(\zeta(x+b)-\zeta(x-b)\big),
\end{equation*}
this will be achieved if we can prove that for $x\rightarrow a$,
\begin{equation*}
\frac{d\log\alpha}{dx}=\frac{1}{x-a}+\frac{w_{-1}}{4}
   -\frac{w_{-1}}{1-{e_4}^2}+O(x-1).
\end{equation*}
This follows by a calculation analogous to the one in Section
\ref{sec:e3=e4}.
\end{proof}
 

\subsection{The Lagrange top as limit case: $e_4=-e_4=\infty$}

Look at the $t$-integral as it is written in equation
(\ref{eqn:tint0}). If one simply lets $e_4$ tend to infinity, the
integrand diverges. But note that from equation
(\ref{eqn:e4}) one obtains 
\begin{equation*}
-\frac{s}{2}=\frac{A}{2p}\,\frac{1}{1-{e_4}^2},
\end{equation*}
so that equation (\ref{eqn:tint0}) can be rewritten as 
\begin{equation*}
t=\sqrt{\frac{A}{2p}}\int\sqrt{\frac{u^2-{e_4}^2}{1-{e_4}^2}}
  \frac{du}{\sqrt{(u-e_1)(u-e_2)(u-e_3)}}.
\end{equation*}
Now if one lets $e_4$ tend to infinity while keeping $A/p$ constant,
this converges to the Abelian integral of the first kind
\begin{equation*}
t=\sqrt{\frac{A}{2p}}\int\frac{du}{w}
\end{equation*}
on the elliptic curve $w^2=(u-e_1)(u-e_2)(u-e_3)$.

The integrals (\ref{eqn:abgdints}) cause no problems. They become
\begin{equation*}
\begin{split}
\log \alpha &= \int \frac{w+w_{-1}}{2(u+1)}\,\frac{du}{w}\\
\log \beta &= \int \frac{w-w_{+1}}{2(u-1)}\,\frac{du}{w}\\
\log \gamma &= \int \frac{w+w_{+1}}{2(u-1)}\,\frac{du}{w}\\
\log \delta&= \int \frac{w-w_{-1}}{2(u+1)}\,\frac{du}{w}\,.
\end{split}
\end{equation*}

This is exactly the solution F. Klein obtains for the Lagrange top
\cite[p. 238]{Klein:Buch},\cite[pp. 28f]{Klein:Lectures}. 
This is not surprising since this limiting case is obtained by letting
$s$ tend to zero while keeping all other parameters
constant. But in this case the Lagrangian (\ref{eqn:lagrange})
coincides with the one of Lagrange's top.


\small

\begin{thebibliography}{10}

\bibitem{IntSysI}
Dubrovin B.A., Krichever I.M. and Novikov S.P.,
\newblock Integrable Systems {I},
\newblock in Dynamical Systems {IV}, Editors V.I. Arnol'd and S.P. Novikov,  in
  Encyclopaedia of Mathematical Sciences, no.~4, Springer, Berlin, 1990,
  173--280.

\bibitem{IntSysII}
Olshanetsky M.A., Perelomov A.M., Reyman A.G. and Semenov-Tian-Shansky M.A.,
\newblock Integrable Systems {II},
\newblock in Dynamical Systems {VII}, Editors V.I. Arnol'd and S.P. Novikov,
  in Encyclopaedia of Mathematical Sciences, no.~16, Springer, Berlin, 1994,
  83--259.

\bibitem{Klein:Lectures}
Klein F.,
\newblock The Mathematical Theory of the Top,
\newblock in Congruence of Sets and other Monographs, Chelsea Publishing
  Company, New York,
\newblock Lectures delivered in Princeton in 1896.

\bibitem{Poisson:Mechanik}
Poisson S.D.,
\newblock Lehrbuch der Mechanik, vol.~2,
\newblock G. Reimer, Berlin, 1836,
\newblock Transl. from 2nd ed.

\bibitem{Whittaker}
Whittaker E.T.,
\newblock A Treatise on the Analytical Dynamics of Particles and Rigid Bodies,
  4th~ed.,
\newblock Cambridge University Press, Cambridge, 1988,
\newblock First published 1904.

\bibitem{Klein:Mechanik}
Klein F.,
\newblock {Einleitung in die analytische Mechanik},
\newblock Teubner, Stuttgart and Leipzig, 1991,
\newblock Lectures, held in G\"ottingen 1886/87.

\bibitem{Klein:Buch}
Klein F. and Sommerfeld A.,
\newblock \"Uber die Theorie des Kreisels,
\newblock B. G. Teubner, Stuttgart, 1965,
\newblock Reprint of the 1897-1910 ed.

\bibitem{BobenkoSuris1}
Bobenko A.I. and Suris Y.B.,
\newblock Discrete Time {L}agrangian Mechanics on {L}ie groups, with an
  Application to the {L}agrange Top,
\newblock {\em Communications in Mathematical Physics}, 1999, V.204, 147--188.

\bibitem{BobenkoSuris2}
Bobenko A.I. and Suris Y.B.,
\newblock Discrete {L}agrangian Reduction, Discrete {E}uler-{P}oincar\'e
  Equations, and Semidirect Products,
\newblock {\em Letters in Mathematical Physics}, 1999, V.49, 79--93.

\bibitem{Klein:Ikosaeder}
Klein F.,
\newblock Vorlesungen \"uber das Ikosaeder und die Aufl\"osung der
  Glei\-chun\-gen vom f\"unften Grade,
\newblock Teubner, Leipzig, 1993,
\newblock Reprint of the 1884 ed.

\end{thebibliography}

\label{lastpage}

\end{document}